# Partial Group Representations on Semialgebras

**Meenakshi, R. P. Sharma**

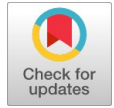

*Abstract: Let A bean additively cancellative semialgebra over an additively cancellative semifield K as defined in [9]. For a given partial action α of a group G on an algebra, the associativity of partial skew group ring together with the existence and uniqueness of enveloping (global) action were studied by M. Dokuchaev and R. Exel [2] which were extended for semialgebras with some restriction by Sharma et. al. using the ring of differences. In a similar way, we extend the results of [2,3] for semialgebras regarding partial representations.*

## I. INTRODUCTION

This paper is in continuation of Sharma et. al. [11], wherein the authors gave some applications of some theorems to partial $(G, \alpha)$-sets and to distinguished the labelling of partial actions. For a given strong partial action $(\alpha, A)$ of a finite group $G$ on an additively cancellative yoked semialgebra $A$, a morita equivalence between Partial Actions and Global Actions between skew group semirings arose from the partial action has been constructed in [10]. Now, exploring different analysis of generalized form of partial actions, Sharma et. al. [9] defined the tensor product of semimodules over semirings and they showed with respect to the defined tensor product the category $(K-Smod, \otimes_K, K)$ becomes a monoidal category, where $K$ is a commutative additively cancellative semiring with identity. The monoids of this monoidal category are $K$-semialgebras. Defining semialgebras, in their paper, the authors studied the global actions of partial actions on semialgebras and the conditions for associativity of partial skew group semirings.

In this paper, we study partial representations of partial actions on semialgebras. We recall from [2] that partial action $\alpha$ of a group $G$ on a semialgebra is defined as:

**Definition 1.1.** A partial action $\alpha$ of $G$ on a set $X$ is a pair $\alpha = \{\{D_g\}_{g \in G}, \{\alpha_g\}_{g \in G}\}$, where for each $g \in G, D_g$ is a subset of $X$ and $\alpha_g : D_{g^{-1}} \to D_g$ is a bijective map, satisfying the following three properties for each $g, h \in G$:

(i) $D_1 = X$ and $\alpha_1 = IdX$, the identity map on $X$,

(ii) $\alpha_g(D_{g^{-1}} \cap D_h) = D_g \cap D_{gh}$,

(iii) $\alpha_g(\alpha_h(x)) = \alpha_{gh}(x)$ for $x \in D_{h^{-1}} \cap D_{h^{-1}g^{-1}}$.

Following [1], for convenience we write $\exists g \cdot x$ to mean that $g \cdot x$ is defined; more precisely, $\exists g \cdot x$ means $x \in D_{g^{-1}}$ and $g \cdot x \in D_g$. Thus from [1], we have an equivalent definition of a partial action of a group on a set as follows:

**Definition 1.2.** The partial action $\alpha$ of a group $G$ on a set $X$ is a partial function from $G \times X$ to $X$ together with the following conditions:

(PA1) $\exists e \cdot x$ for all $x \in X$ and $e \cdot x = x$, where $e$ is the identity of $G$;

(PA2) $\exists g \cdot x$ implies that $\exists g^{-1} \cdot (g \cdot x)$ and $g^{-1} \cdot (g \cdot x) = x$, $g \in G$ and $x \in X$;

(PA3) $\exists g \cdot (h \cdot x)$ implies that $\exists gh \cdot x$ and $g \cdot (h \cdot x) = (gh) \cdot x$.

In case $X$ is a semialgebra, partial maps are semialgebra homomorphisms. Throughout this paper, $G$ is a finite group acting partially on an additively cancellative $K$-semialgebra $A$, where $K$ is an additively cancellative semifield. The identities of both $A$ and $K$ are denoted by 1 and the identity of $G$ by $e$.

If $G$ is a finite group and $K$ is a field, then representations theory of $G$ on $K$- algebras is identical to the representation of associative group algebra $KG$. The associative partial group algebra of $G$, denoted by $Kpar(G)$ plays similar role for the partial representations of $G$. In particular, for a finite group $G$, there is a one-one correspondence between the partial representations of $G$ and the representations of $Kpar(G)$, which is isomorphic to groupoid semialgebra $K\Gamma(G)$. These results are proved in section 3.

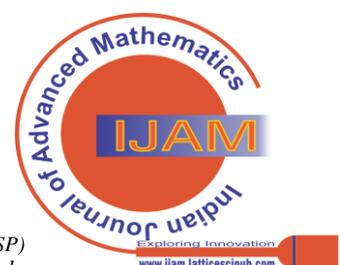

In section 4, we show that $K\Gamma(G)$ is a direct sum of matrix semialgebras over the semirings $KH$, where $H$ is a subgroup of $G$.

## II. PRELIMINARIES

First we recall the definition of a variety from Mac Lane S. [7][8].

**Definition 2.1.** (Variety of algebras) An algebraic system of a type $\tau$ is given by a set $\Omega$ of operators and a set $E$ of identities. The set $\Omega$ of operators is a graded set; that is, a set $\Omega$ with a function which assigns to each $\omega \in \Omega$ a non negative integer $n$ called the arity of $\omega$. From the given operators $\Omega$, we can form the set $\Lambda$ of all derived operators; given $\omega$ of arity $n$ and $n$ derived operators $\lambda_1, \lambda_2, \ldots, \lambda_n$ of arities $m_1, m_2, \ldots, m_n$, the "composite" $(\lambda_1, \lambda_2, \ldots, \lambda_n)\omega$ is a derived operator of arity $m_1 + m_2 + \ldots + m_n$; also given $\lambda$ of arity $n$ and $f : n \to m$ any function from $\{1, 2, \ldots, n\}$ to $\{1, 2, \ldots, m\}$, substitution of $f$ in $\lambda$ gives a derived operator $\theta$ of arity $m$, described in terms of variables $x_i$ as $\theta(x_1, x_2, \ldots, x_m) = \lambda(f_{x_1}, f_{x_2}, \ldots, f_{x_n})$. The set $E$ of identities for algebraic systems of type $\tau$ is a set of ordered pairs $\langle \lambda, \mu \rangle$ of derived operators, where $\lambda$ and $\mu$ have the same arity $n$. If $S$ is any set, an action of $\Omega$ on $S$ is a function $f$ which assigns to each operator $\omega$ of arity $n$ an $n$−ary operation $\omega_f : S^n \to S$.

An action $f$ of $\Omega$ on $S$ satisfies the identity $\langle \lambda, \mu \rangle$ if $\lambda_f = \mu_f : S^n \to S$.

An $\langle \Omega, E \rangle$−algebra is a set $S$ together with an action $f$ of $\Omega$ on $S$ which satisfies all the identities of $E$.

A morphism $g : f \to f'$ of $\langle \Omega, E \rangle$−algebras is a function $g : S \to S'$ on the underlying sets which preserves all the operators of $\Omega$ in the sense that $g((a_1, a_2, \ldots, a_m) \omega_f) = (ga_1, ga_2, \ldots, ga_m)\omega_{f'}$ for all $a_i \in S$, where $f$ and $f'$ are actions of $\Omega$ on $S$ and $S'$ respectively.

The collection of all small $\langle \Omega, E \rangle$−algebras, with their morphisms as arrows, is a category $\langle \Omega, E \rangle - Alg$, often called a variety.

**Definition 2.2.** (Katsov Y.[6]) An algebra $A$ of type $\tau$ is called abelian if for each algebra $B$ of type $\tau$ the following conditions hold

(i) For any $n - ary$ ($n \geq 1$) operator $\omega$ and morphisms of algebras $\phi_1, \phi_2, \ldots, \phi_n$ from $B$ to $A$, the function $\phi_1 \phi_2 \ldots \phi_n \omega$ defined by
$(\phi_1 \phi_2 \ldots \phi_n \omega)b = (\phi_1 b)(\phi_2 b) \ldots (\phi_n b) \omega$ for $b \in B$
is a morphism from $B$ to $A$.

(ii) For any nullary operator $\omega$, the function $\phi\omega : B \to A$ defined by $(\phi\omega) b = 0\omega, b \in B$, where $0\omega$ is the element in $A$ distinguished by $\omega$; is a morphism of algebras.

A variety whose all algebras are abelian is called an abelian variety.

**Example 2.1.** Let $\widetilde{M}$ the category of all commutative monoids. Then $\widetilde{M}$ is a variety. Here $\Omega$ contains two operators, the product and the assignment of identity element $e$, of arities 2 and 0 respectively, and $E$ contains (i) the axiom for the identity ($ex = x = xe$); (ii) the associative and commutative laws.

Let $F, G \in \widetilde{M}$ and $F(F, G)$ denote the free monoids generated by $F \times G$, $\rho$ be the congruence on $F(F, G)$ generated by the pairs $\langle (a_1 a_2, b), (a_1, b)(a_2, b) \rangle$ and $\langle (a, b_1 b_2), (a, b_1)(a, b_2) \rangle$. Take $F \otimes G$ as $F(F, G)/\rho$. Let $F_1, G_1 \in \widetilde{M}$ and $\alpha \in \widetilde{M}(F, F_1)$ and $\beta \in \widetilde{M}(G, G_1)$, the assignments $(F, G) \mapsto F \otimes G$ and $(\alpha, \beta) \mapsto \alpha \otimes \beta$ determine the bifunctor $\_ \otimes \_ : \widetilde{M} \times \widetilde{M} \to \widetilde{M}$. Since $\widetilde{M}$ is an abelian variety. So, the bifunctor $\_ \otimes \_ : \widetilde{M} \times \widetilde{M} \to \widetilde{M}$ is an internal tensor product by Katsov Y. [6].

Now we recall some results and definitions from [9].

Let $K - Smod$ and $Smod - K$ denote the categories of left and right $K-$ semi- modules, respectively, over a semiring $K$. The following three results can be derived easily in the category $\widetilde{M}$ of commutative monoids in place of commutative inverse monoids and without assuming $K$ an additively regular semiring as in Katsov Y.

**Proposition 2.1.** Let N be a commutative monoid, $G_1, G \in K - Smod$ and $\eta \in K - Smod(G_1, G$. Let $\widetilde{M}(\eta, 1_N) : \widetilde{M}(G, N) \to \widetilde{M}(G_1, N)$ be the map

$\phi \to 1_N \circ \phi \circ \eta$ for $\phi \in \widetilde{M}(G, N)$. Then the assignments $G \to \widetilde{M}(G, N)$ and $\eta \to \widetilde{M}(\eta, 1_N)$

determine the additive functor $N$ from $(K - Smod)^{op}$ to $Smod - K$.

**Proposition 2.2.** Let G be a left $K-$semimodule, $N, M \in \widetilde{M}$ and $f \in \widetilde{M}(N, M)$. Let $\widetilde{M}(1_G, f) : \widetilde{M}(G, N) \to \widetilde{M}(G, M)$ be the map $\theta \to f \circ \theta \circ 1_G$ for $\theta \in \widetilde{M}(G, N)$. Then the assignments $N \mapsto \widetilde{M}(G, N)$ and $f \mapsto \widetilde{M}(1_G, f)$ determine the additive functor $(\_)^G$ from $\widetilde{M}$ to $Smod - K$.

**Proposition 2.3.** For $G, G_1 \in K - Smod, N, M \in \widetilde{M}, \eta \in K - Smod(G_1, G)$ and $f \in \widetilde{M}(N, M)$, let $\widetilde{M}(\eta, f) : \widetilde{M}(G, N) \to \widetilde{M}(G_1, M)$ be the map $\theta \mapsto f \circ \theta \circ 1_G$ for $\theta \in \widetilde{M}(G, N)$.

21

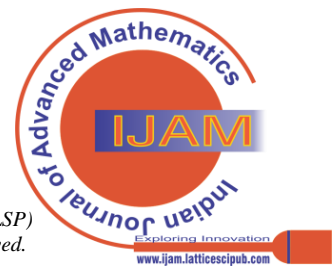

Then the assignments $(G, N) \mapsto \widetilde{M}(G, N)$ and $(\eta, f) \mapsto \widetilde{M}(\eta, f)$ determine a bifunctor from $(K - Smod)^{op} \times \widetilde{M}$ into $Smod - K$.

**Definition 2.3.** (Tensor product of semimodules) Let $F \in Smod - K$ and $G \in K - Smod$ then both $F$ and $G$ are commutative monoids, so are in $\widetilde{M}$ and therefore has tensor product $F \otimes G$ (considered as commutative monoids). The tensor product $F \otimes_K G$ is defined as the factor monoid $(F \otimes G)/\sigma$, where $\sigma$ is the congruence on $F \otimes G$ generated by the pairs $\langle xk \otimes y, x \otimes ky \rangle$ for all $x \in F, y \in G$ and $k \in K$.

Given $F, G \in \widetilde{M}$, let $bih(F, G)$ denote the category having bihomomorphisms $f : F \times G \to C$ as objects; and morphisms between bihomomorphisms $f : F \times G \to C$ and $f_1 : F \times G \to C_1$ are homomorphisms $\alpha : C \to C_1$ of $\widetilde{M}$ such that $\alpha \circ f = f_1$.

**Lemma 2.1.** If $F \in Smod - K$ and $G \in K - Smod$, then for any balanced product $(C, f)$ of $F$ and $G$ there exists a unique morphism of monoids $\phi : F \otimes_K G \to C$ such that $f = \phi \circ g$, where $g$ is the composition of the morphism $F \times G \xrightarrow{h} F \otimes G \xrightarrow{\pi} F \otimes_K G$; $h$ is the initial object in the category $bih(F, G)$ and $\pi$ is the canonical epimorphism.

The usual tensor product [5] is not valid for semimodules over semirings and the category $K - Smod$ is not a monoidal category with the respect to the tensor product defined in [4]. However, the tensor product defined in [10] gives

**Theorem 2.1.** Let $K$ be a commutative semiring. Then $(K - Smod, \otimes K, K)$ is a monoidal category.

**Definition 2.4.** The monoids in the monoidal category $(K - Smod, \otimes K, K, \alpha, \lambda, \varrho)$ are called $K$−semialgebras.

### III. PARTIAL GROUPOID SEMIALGEBRAS

Let $G$ be a finite group with the identity $e$ and $K$−an additively cancellative semi-field. The partial groupoid semialgebra $K_{par}(G)$ is defined similar to the partialgroupoid algebra. In this section we consider groupoid $\Gamma(G)$ constructed for $G$ in[3] which gives $K\Gamma(G) \cong K_{par}(G)$ and there is a one-one correspondence between the partial representations of $G$ and the representations of $K\Gamma(G)$. First, we define

**Definition 3.1.** A partial representation of a group $G$ into a unital $K$−semialgebra $A$ is a map $\pi : G \to A$ such that for all $g, h \in G$, we have $\pi(g)\pi(h)\pi(h^{-1}) = \pi(gh)\pi(h^{-1}), \pi(g^{-1})\pi(g)\pi(h) = \pi(g^{-1})\pi(gh)$ and $\pi(e) = 1$, where $e$ is the identity of $G$. Denote by $K_{par}(G)$, the universal $K$−semialgebra generated by symbols $\{[g], g \in G\}$ satisfying the following conditions:

(1) $[e] = 1$;

(2) $[s^{-1}][s][t] = [s^{-1}][st]$;

(3) $[s][t][t^{-1}] = [st][t^{-1}]$; for all $s, t \in G$.

Remark 3.1. The map $g \mapsto [g]$ is a partial representation of on $K_{par}(G)$. If A is an additively cancellative $K$−semialgebra and $\pi : G \to A$ is a partial representation of G on A, then π extends uniquely by linearity to a representation $\phi: Kpar(G) \to A(\phi([g]) = \pi(g))$. Conversely, if $\phi : Kpar(G) \to A$ is a unital homomorphism of $K$−semialgebra, then $\pi(t) = \phi([t])$ is a partial representation of G on A.

We show that for a groupoid $\Gamma(G)$, there is a one-one correspondence between the partial representations of G and the representations of groupoid semialgebra KΓ(G). First we recall the definition of $\Gamma(G)$. from [3]. For a finite group $G$, $\Gamma(G) = \{(I, g) | g \in G$ and $e, g^{-1} \in I \subseteq G\}$, so for $(I, g) \in \Gamma(G)$ we have e, $g \in gI$. The multiplication of pairs $(I, g), (J, h)$ in $\Gamma$ is defined if $I = hJ$ and in this case

$$(hJ, g) \cdot (J, h) = (J, gh)$$

The groupoid $K$—semialgebra KΓ(G) is given by $\sum_i \alpha_i \gamma_i, \alpha_i \in K, \gamma_i \in \Gamma\}$ in which the sum is component wise and the multiplication is given by:

For $\alpha(I, g), \beta(J, h) \in K\Gamma(G)$,

$$\alpha(I, g) \cdot \beta(J, h) = \begin{cases} \alpha\beta(J, gh), \text{if } I = hJ \\ 0, \text{ otherwise} \end{cases}.$$

Note that if $K$ is additively cancellative, then $K\Gamma(G)$ is additively cancellative. Hence, both $K^\triangle$ and $(K\Gamma(G))^\triangle$ exist and in this case we have

**Lemma 3.1.** $K^\triangle\Gamma(G)$ is a groupoid semialgebra and $(K\Gamma(G))^\triangle \cong K^\triangle\Gamma(G)$.

**Proof.** Since $K\Gamma(G)$ is $K$−semialgebra, by [9] $(K\Gamma(G))^\triangle \cong K\Gamma(G)^\Delta$ is $K^\Delta$−algebra. Hence, it suffices to show that $(K\Gamma(G))^\triangle \cong K^\triangle\Gamma(G)$.

Define a map $f: K^\Delta\Gamma(G) \cong K\Gamma(G)^\Delta$ by

$$f\sum_i((\alpha_i - \beta_i)\gamma_i) = \sum_i \alpha_i\gamma_i - \sum_i \beta_i\gamma_i \quad (3.1)$$

Since $\{\gamma_i\}$ are $K^\Delta$−basis for $K^\Delta\Gamma(G)$ and $K$−basis for $K\Gamma(G), \sum_i((\alpha_i - \beta_i)\gamma_i = \sum_i(\alpha'_i - \beta_i')\gamma_i)$, implies $(\alpha_i - \beta_i) = (\alpha'_i - \beta'_i)$ for all $i$, so that $(\alpha_i + \beta_i') = (\alpha'_i + \beta_i)$. Therefore, $(\alpha_i\gamma_i + \beta'_i\gamma_i) = (\alpha_i'\gamma_i + \beta_i\gamma_i)$ and hence $\sum_i \alpha_i\gamma_i - \sum_i \beta_i\gamma_i = \sum_i \alpha'_i\gamma_i - \sum_i \beta_i'\gamma_i$, proving that the map is well-defined. By reversing the steps, it follows that $f$ is one-one. For onto, let $\sum_i \alpha_i \gamma_i - \sum_i \beta_i \gamma_i \in (K\Gamma(G))^\Delta$. Then $\sum_i ((\alpha_i - \beta_i)\gamma_i) \in K^\Delta\Gamma(G)$ and

$$f(\sum_i((\alpha_i - \beta_i)\gamma_i) = \sum_i \alpha_i \gamma_i - \sum_i \beta_i \gamma_i$$

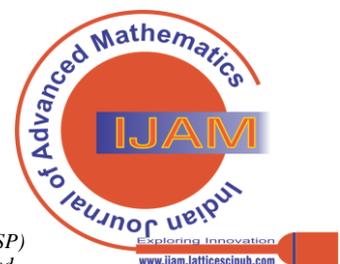

Now, we will show that $f$ is a homomorphism. For, let $\gamma_1, \gamma_2 \in \Gamma(G)$, where $\gamma_1 = (hJ, g)$ and $\gamma_2 = (J, h)$, then

$$f((\alpha_1 - \beta_1)\gamma_1 \cdot (\alpha_2 - \beta_2)\gamma_2) = f((\alpha_1 - \beta_1)(hJ,g) \cdot (\alpha_2 - \beta_2)(J,h)) = f((\alpha_1 - \beta_1)(\alpha_2 - \beta_2)(hJ,g)(J,h)) = f((\alpha_1 - \beta_1)(\alpha_2 - \beta_2)(J,gh)) = f((\alpha_1\alpha_2 + \beta_1\beta_2) - (\alpha_1\beta_2 + \beta_1\alpha_2))(J,gh)) = (\alpha_1\alpha_2 + \beta_1\beta_2)(J,gh) - (\alpha_1\beta_2 + \beta_1\alpha_2)(J,gh) = (\alpha_1 - \beta_1)(hJ,g) \cdot (\alpha_2 - \beta_2)(J,h) = f((\alpha_1 - \beta_1)\gamma_1)f(((\alpha_2 - \beta_2)\gamma_2).$$

We note from the above lemma that $K\Gamma(G)$ embeds in $K^\Delta\Gamma(G)$ as $K$ embeds $K^\Delta$ and the elements $(I, e) \in \Gamma(G)$ are mutually orthogonal and their sum is the identity in both $K\Gamma(G)$ and $K^\Delta\Gamma(G)$.

**LeΣmma 3.2.** The map $\lambda p: G \to K\Gamma(G \to K^\Delta\Gamma(G)$ defined by $\lambda p(g) = \sum_{g^{-1}\in I}(I,g)$ is a partial representation of $G$.

**Proof.** This follows by using the multiplication in $\Gamma(G)$ and the fact that $\sum_{e\in I}(I,e) = 1$.

**Theorem 3.1.** There is a one-to-one correspondence between the partial representations of G and the representations of $K\Gamma(G)$. More precisely, if A is any additively cancellative unital K-semialgebra, then $\pi : G \to A$ is a partial representation of G if and only if there exists a unital semialgebra homomorphism $\tilde{\pi}: K\Gamma(G) \to A$ such that $\pi = \tilde{\pi} \circ \lambda_p$. Moreover, such a homomorphism $\tilde{\pi}$ is unique.

**Proof.** If $\tilde{\pi}: K\Gamma(G) \to A$ is a unital homomorphism of $K -$ semialgebras, then we have to show that the composite map $\pi = \tilde{\pi}\circ\lambda p: G \to A$ is a partial representation of $G$ on $A$. Since $\tilde{\pi}$ is a unital semialgebra homomorphism $\tilde{\pi}: K\Gamma(G) \to A$, therefore $\pi(e) = 1$. For $g, h \in G$, we have $\pi(g^{-1})\pi(g)\pi(h) = \tilde{\pi}(\sum_{g\in I}(I,g^{-1}))(\tilde{\pi}(\sum_{g^{-1}\in I}(I,g)\tilde{\pi}\sum_{h^{-1}\in J}(J,h)) = \tilde{\pi}(\sum_{g^{-1}\in I}(I,e)\tilde{\pi}\sum_{h^{-1}\in J}(J,h)) = \tilde{\pi}(\sum_{g^{-1}\in I}(I,g)\sum_{h^{-1}\in J}(J,h)) = \tilde{\pi}(\sum_{h^{-1}g^{-1}\in J}(J,h)).$

On the other hand, $\pi(g^{-1})\pi(gh) = \tilde{\pi}(\sum_{g\in I}(I,g^{-1}))(\tilde{\pi}(\sum_{h^{-1}g^{-1}\in D}(D,gh)) = \tilde{\pi}(\sum_{h^{-1}g^{-1}\in D}(ghD,g^{-1})(D,gh)) = \tilde{\pi}(\sum_{h^{-1}g^{-1}\in D}(D,h))$, so that $(g^{-1})\pi(g)\pi(h) = \pi(g^{-1})\pi(gh)$. Similarly, we can verify that $\pi(g)\pi(h)\pi(h^{-1}) = \pi(gh)\pi(h^{-1})$.

Conversely, suppose $\pi: G \to A \to A^\Delta$ be a partial representation of $G$ on $A$, so by [3, Theorem 2.6], there exists a unital homomorphism $\tilde{\pi}^\Delta: K^\Delta\Gamma(G) \to A^\Delta$ such that $\pi = \tilde{\pi}^\Delta \circ \lambda_p$.

We note from [3, Theorem 2.6] that for all elements $(I, g)$ in $\Gamma(G)$, we have $\tilde{\pi}^\Delta(I,g) = \pi(g) \prod_{r\in I}\epsilon(r) \prod_{s\in I}(1 - \epsilon(s))$, where $\epsilon(r) = \pi(r)\pi(r^{-1}), g \in G$. Hence, $\tilde{\pi}^\Delta(I,g) \in A$, since $\pi$ is initially a map from $G$ to $A$. As $K^\Delta\Gamma(G) \cong (K\Gamma(G))^\Delta$, let $\tilde{\pi}$ be the restriction of $\tilde{\pi}^\Delta$ on $K\Gamma(G)$. Then $\tilde{\pi}(\sum \alpha_i(I,g_i)) = \sum \alpha_i \tilde{\pi}(I,g_i) \in A$, as $\tilde{\pi}^\Delta(I,g) \in A$ for all $(I,g) \in \Gamma G)$.

Hence, we have the map $\tilde{\pi}^\Delta: K\Gamma(G) \to A$, a homomorphism which obviously satisfies $\tilde{\pi}^\Delta(\lambda_p(g)) = \pi(g)$. Suppose there exists an another homomorphism $\tau: K\Gamma(G) \to A$ such that $\pi = \tau \circ \lambda p$ with $\tilde{\pi}' \neq \tilde{\pi}$ Then $\tau$ can be extended to a homomorphism $\tau^\Delta((\alpha - \beta)\gamma) = \tau(\alpha\gamma) - \tau(\beta\gamma), \alpha, \beta \in K, \gamma \in \Gamma(G)$. Obviously, $\tau^\Delta$ also satisfies $\pi = \tau^\Delta \circ \lambda_p$. This contradicts the uniqueness of $\tilde{\pi}^\Delta$ as $\tau \neq \tilde{\pi}$ implies $\tau^\Delta \neq \tilde{\pi}^\Delta$.

**Corollary 3.1.** The groupoid semialgebra $K\Gamma(G)$ is isomorphic to the partial group semialgebra $K_{par}(G)$.

**Proof.** The maps $[]: G \to K_{par}(G)$, given by $g \to [g]$ and $\lambda p : G \to K\Gamma(G)$ are partial representation of $G$. By the universal property of the two semialgebras $K_{par}(G)$ and $K\Gamma(G)$, there exists a unital $K -$semialgebra homomorphism $\tilde{\pi}: K\Gamma(G) \to K_{par}(G)$ and $\phi : K_{par}(G) \to K\Gamma(G)$ such that $\tilde{\pi}(\lambda p(g)) = [g]$ and $\phi([g]) = \lambda_p(g)$.

Now, for the composite map $\tilde{\pi} \circ \phi : K_{par}(G) \to K_{par}(G)$, we have $\tilde{\pi}\circ\phi([g]) = \tilde{\pi}(\lambda_p(g)) = [g]$

And for the composite map $\phi \circ \tilde{\pi}: K\Gamma(G) \to K\Gamma(G)$, we have $\phi \circ \tilde{\pi}(\sum_{g^{-1}\in I}(I,g)) = \phi \circ \tilde{\pi}(\lambda_p(g)) = \phi(\pi(g)) = \lambda_p(g) = \sum_{g^{-1}\in I}(I,g)$, that is, these maps are the identities on their domains, since $\{[g] : g \in G\}$ is the set of generators for $K_{par}(G)$ and $\{(I,g) \in \Gamma(G)\}$ is the set of generators for $K\Gamma(G)$.

## IV. STRUCTURE OF THE PARTIAL GROUP SEMIALGEBRAS

In this section, we study the structure of the groupoid semialgebra $K\Gamma(G)$ which is isomorphic to $K_{par}(G)$; we prove that $K\Gamma(G)$ is the direct sum of matrix semialgebra over the semirings $KH$, where $H$ is any subgroup of $G$. We follow the following notation from [3].

Given a finite group $H$ and a positive integer $m$, let $\Gamma_m^H$ denote the groupoid whose elements are triplets $(h, i, j)$ where $h \in H$ and $i, j \in (1, 2, \ldots, m)$. The source and the range map on $\Gamma_m^H$ are defined by $s(h,i,j) = j$ and $r(h,i,j) = i$. The product is defined by $(g,i,j) \cdot (h,j,k) = (gh,i,k)$. The units of $\Gamma_m^H$ are the elements of the form $(e,i,i); i = 1,2,\ldots,m$.

A groupoid $\Gamma$ is represented with an oriented graph $E_\Gamma$, whose vertices are the units of the groupoid, and each element $\gamma \in \Gamma$ gives an oriented edge of $E_\Gamma$ from the vertex $s(\gamma)$ to the vertex $r(\gamma)$.

23

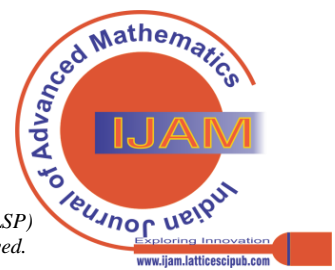

A connected component of $E_\Gamma$ gives a subgroupoid of $\Gamma$.

In the case of the groupoid $\Gamma_m^H$, the corresponding graph $E_{\Gamma_m^H}$ has precisely $m$ vertices, and between any two vertices there are precisely $|H|$ oriented edges (in each direction) labeled by the elements of H.

**Proposition 4.1.** Let $\Gamma$ be a groupoid such $E_\Gamma$ that is connected and $m = |\Gamma^{(0)}|$ is finite. Let $x_1$ be any vertex of $E_\Gamma$, and $H$ be the isotropy group of $x_1$, which is defined by $H = \{\gamma \in \Gamma : s(\gamma) = r(\gamma) = x_1\}$. Then

$$(1) \quad \Gamma \cong \Gamma_m^H$$

$$(2) \quad K\Gamma \cong M_m(KH).$$

The proof of (1) is same as in [3].

Before proving (2), we prove a result which will be used to define an isomorphism between groupoid semialgebra $K\Gamma$ and $M_m(KH)$.

**Lemma 4.1.** Let $K$ be an additively cancellative semifield and $H$ be any groupoid. Then $M_m(KH) \cong M_m(K) \otimes_K KH$, where $\otimes_K$ is defined as in [9].

**Proof.** Since $M_m(K)$ is a right $K$-semimodule over $K$ and $KH$ is a left $K$-semimodule over $K$, we define a map $f: M_m(K) \times KH \to M_m(KH)$, where $M_m(KH)$ is an additive abelian group, such that $f(A, \alpha h) = (\alpha_{i,j}\alpha h)_{m \times n}$.

Then $(M_m(KH), f)$ is a balanced product of $M_m(K)$ and $KH$. Therefore, by Lemma 2.1, there exists a unique homomorphism of semimodules $\phi: M_m(K) \otimes_K KH \to M_m(KH)$ defined by $\phi(A \otimes_K h) = f(A, h)$ such that $f = \phi \circ i$, where $i$ is the map from $M_m(K) \times KH$ to $M_m(K) \otimes_K KH$ defined by $(a, b) \mapsto (a \otimes b)$.

Now, we have to prove that $\phi$ is one-one and onto. For, define a map $\varphi: M_m(KH) \to M_m(K) \otimes_K KH$ given by $\varphi([a_{ij}]_{m \times m}) = \sum_{i,j}(e_{ij} \otimes a_{ij})$ where $\{e_{ij}), i,j = 1, 2, \ldots, m\}$ denotes the set of matrix units of $M_m(K$;i.e., $e_{ij}$ is the matrix with null entries, except for the entry in the position $(i, j)$ which is equal to the identity of $K$. Clearly, $\varphi$ is a homomorphism of $K$-semimodules.

Now,

$$\phi \circ \varphi([a_{ij}]_{m \times m}) = \phi(\varphi([a_{ij}]_{m \times m})$$
$$= \phi\left(\sum_{i,j}(e_{ij} \otimes a_{ij})\right)$$
$$= [a_{ij}]_{m \times m}$$

And $\varphi \circ \phi = \varphi(\phi(A \otimes \sum_i \alpha_i h_i)) = \varphi(f(A, \sum_i \alpha_i h_i)) = \varphi([b_{ij}]_{m \times m})$, where $b_{ij} = a_{ij} \sum_i \alpha_i h_i = \sum_{i,j}(e_{ij} \otimes b_{ij}) = A \otimes \sum_i \alpha_i h_i$, that is the two composite maps $\phi \circ \varphi$ and $\varphi \circ \phi$ are the identity maps on their respective domains. Therefore, $\varphi$ is the inverse of $\phi$; hence $\phi$ is one-one and onto, from which we have $Mm(KH) \cong Mm(K) \otimes_K KH$.

**Proof of Lemma 4.1.(2).** Let $\{(eij), i, j = 1, 2, \ldots, m\}$ be as in above lemma. For any group $H$, the map $(h, i, j) \mapsto e_{i,j} \otimes h$ extends by linearity to a $K$-semialgebra isomorphism between $K\Gamma_m^H$ and $Mm(K) \otimes_K KH \cong Mm(KH)$, which gives the required result.

**Theorem 4.1.** The semigroupoid semialgebra $K\Gamma(G)$ is of the form $K\Gamma(G) = \bigoplus_{\substack{H \le G \\ 1 \le m \le (G:H)}} c_m(H)M_m((KH)$

($c_m(H)M_m((KH)$ means the direct sum of $c_m(H)$ copies of $M_m((KH)$, where

$$C_m(H) = \left(\binom{(G:H) - 1}{m - 1}\right) - m \sum_{\substack{H < B \le G \\ (B:H)|m}} \frac{(c_{\frac{m}{B;H}}(B)}{(B:H)}$$

**Proof.** Let $I$ be any subset of $G$ containing the identity of $G$. Let $H = S(I)$, the stabilizer of $I$ in $G$, given by $S(I) = \{g \in G: gI = I\}$. In the graph $E_{\Gamma(G)}, S(I)$ is identified with the set of edges departing and terminating at the vertex $(I, e)$. Since $e \in I$, we have $H \subseteq I$. Since $H$ acts on the left on $I$, then $I$ is the union of right cosets of $H$, say $I = \bigcup_i^m Ht_i$, $t_1 = e$, where $m = \frac{|I|}{|H|}$. As observed in [3, Theorem 3.2] the sub-groupoid of $\Gamma(G)$ corresponding to the connected component of the vertex $I$ of the graph $E_{\Gamma(G)}$ is isomorphic to the groupoid $\Gamma_m^H$. Observe that the stabilizer of $I$, which is $H$, coincides with isotropy group of the unit $(I, e) \in \Gamma(G)^{(0)}$. Hence, by Proposition 4.1, the semialgebra $M_m(KH \simeq K\Gamma_m^H$ is a direct summand of $K\Gamma(G)$ and

$$K\Gamma(G) = \bigoplus_{\substack{H \le G \\ 1 \le m \le (G:H)}} c_m(H)M_m((KH).$$

$$\left(\bigoplus_{\substack{H \le G \\ 1 \le m \le (G:H)}} c_m(H)M_m((KH)\right)^\Delta, \text{which gives}$$

$$K^\Delta \Gamma(G) = \bigoplus_{\substack{H \le G \\ 1 \le m \le (G:H)}} c_m(H)(M_m((KH))^\Delta \cong$$

$$\bigoplus_{\substack{H \le G \\ 1 \le m \le (G:H)}} c_m(H)(M_m((K^\Delta H)).$$

The rest of the proof follows from equation (14) of [3].

## V. CONCLUSION

For a finite group $G$, there is a one-one correspondence between the partial representations of $G$ and the representations of $Kpar(G)$, which is isomorphic to groupoid semialgebra $K\Gamma(G)$. We see that $K\Gamma(G)$ is the direct sum of matrix semialgebras over the semirings $KH$, where $H$ is a subgroup of $G$.

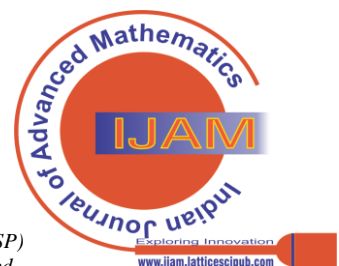

# Partial Group Representations on Semialgebras

## DECLARATION

| | |
|---|---|
| Funding/ Grants/ Financial Support | No, I did not receive. |
| Conflicts of Interest/ Competing Interests | No conflicts of interest to the best of our knowledge. |
| Ethical Approval and Consent to Participate | The article does not require ethical approval and consent to participate with evidence. |
| Availability of Data and Material/ Data Access Statement | Not relevant. |
| Authors Contributions | All authors have equal participation in this article |

## REFERENCE

1. Keunbae Choi and Lim Yongdo, Transitive partial actions of groups, *Periodica Mathematica Hungarica*, Vol 56 (2), 2008, pp. 169-181. [CrossRef]
2. M. Dokuchaev and R. Exel, Associativity of crossed products by partial actions, enveloping actions and partial representations, *Trans. Amer. Math.* Soc. 357(5) (2005) 1931-1952. [CrossRef]
3. M. Dokuchaev, R. Exel and P. Piccione, Partial Representations and Partial group algebras, *Journal of Algebra* 226, 505-532(2000). [CrossRef]
4. J. S. Golan, *Semirings and their applications* (Kluwer Academic Publishers, 1999). [CrossRef]
5. N. Jacobson, *Basic Algebra II*, Hindustan Publishing Corporation (India).
6. Y. Katsov (1997), Tensor Products and Injective Envelopes of Semimodules over Additively Regular Semirings. *Algebra Colloquium* 4:2 (1997) 121-131, 1997.
7. S. Mac Lane, *Categories for the working Mathematician,* Springer-Verlag, New York, 1971.
8. R.P. Sharma and Anu, Semialgebras and their algebras of differences with partial group actions on them, *Asian-Eur. J. Math.,* vol. 6, no. 3 (2013) 1350038 (20 pages). [CrossRef]
9. R.P. Sharma, Anu and N. Singh, Partial group actions on semialgebras, *Asian- Eur. J. Math.* 5(4) (2012), Article ID:1250060, 20pp. [CrossRef]
10. R.P. Sharma and Meenakshi, Morita Equivalence between Partial Actions and Global Actions for Semialgebras, J. of Combinatrics, Information & System Sciences, 41(1-2), 2016, 65-77.
11. R.P. Sharma, Rajni Parmar and Meenakshi, Distinguishing Labelling of Partial Actions, International Journal of Algebra, (9)8, 2015, 371-377. [CrossRef]

## AUTHORS PROFILE

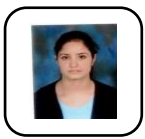

**Dr. Meenakshi,**
Designation: Assistant Professor
Education: Ph. D.
Area of Interest: Algebra
Teaching Experience: 2 years
Number of supervised M.Sc. students for their Project Submission: 12

- R. P. Sharma, Rajni Parmar and Meenakshi, Distinguishing Labelling of Partial Actions, International Journal of Algebra, 9 (8) (2015), 371-377.
- R. P. Sharma and Meenakshi, Morita Equivalence Between Partial Actions and Global Actions for Semialgebras, J. of Combinatorics, Information & System Sciences, 41 (No.1-2) (2016), 65-77.
- R. P. Sharma and Meenakshi, On Construction of Global Actions for Partial Actions, to appear in Proceedings of ICAA to be published by De Gruyter, 2019, pn 243-251.
- Meenakshi, Disjoint Union Metric Spaces, J. of Combinatorics, Information & System Sciences, 2017, 42, pn 81-87.
- R. P. Sharma, Meenakshi and Natakshi, Disjoint Union Metric and Topological Spaces, Southeast Bulletin of Mathematics, 2020, ,44 pg 733-753.
Conference Proceedings: 03

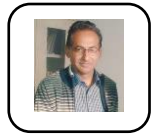

**Prof. R. P. Sharma**
Education: Ph.D
Area if Interest: Algebra
Experience: More than 30 years
Number of supervised Ph.D. students: 12
Some Significant Publications:

1. M. Parvathi and R.P. Sharma, Prime Kernel Functors of group graded rings and their Smash Products, Comm. in Algebra (USA), 17(7) (1989), 1533-1562.
2. M. Parvathi and R.P. Sharma, A note on algebraic analogs of the cones spectrum, Comm. in Algebra (USA), 19(4) (1991), 1249-1270.
3. R.P. Sharma and Samriti Sharma, Group Action on Fuzzy Ideals, Comm. in Algebra (USA), 26(12) (1998) 4207-4220. 4. R.P. Sharma and Samriti Sharma, Graded Fuzzy Ideals in Graded Rings, Ganita, Vol. 49, No. 1 (1998) 29-37.
5. R.P. Sharma and Samriti Sharma, On the Range of A G-Prime Fuzzy Ideals, Comm. in Algebra (USA), 27(6) (1999) 2914-2916.
6. R.P. Sharma, J.R. Gupta and Arvind, Characterization of G-Prime Fuzzy ideals in a ring - An alternate approach, Comm. in Algebra (USA), 28(10) (2000) 1981-1987.
7. R.P. Sharma and Samriti Sharma, Weakly Graded-Prime and Weakly G-Prime Fuzzy ideals, EJM (2000) 1-9.
8. R.P. Sharma and Samriti Sharma, On the range of a Graded-Prime Fuzzy ideal, Ganita, 51, No. I (2000) 69-74. 9. R.P. Sharma, Basic Concepts of Ring Theory (Expository Article), Proceedings of the Institutional and Instructional Program on Quantum Groups and Their Applications, University of Madras, Chennai, 01-02 ( 2001) 39-57.
10. R.P. Sharma, Torsion theoretic Connes subgroup, Ganita, Vol. 52, No2, (2001) 127-135.
11. R.P. Sharma and Ranju Banota, Fuzzy Ideals of Group Graded Rings and their Smash Products, Comm. in algebra (USA), 30(6) (2002)..
12. R.P. Sharma, On the Semi simplicity of Brauer Algebras, University of Glasgow, Department of Mathematics, Preprint Series, Paper No. 2002/45, September 2002.
13. R.P. Sharma, J.R. Gupta and Arvind, On the Characterization of Completely Prime and Completely G-Prime Fuzzy Ideals, Chapter 7, A Book Edited by R.C. Sharma, Published by Allied Publishers Pvt. Ltd. New Delhi (2003) 87-100.